\documentclass[reqno, psamsfonts, 11pt]{amsart}
\usepackage{mathabx}    
\usepackage{mathrsfs}  
\usepackage[backref,colorlinks]{hyperref} 
\usepackage{amsmath,amssymb,amsthm}
\usepackage[abbrev,msc-links]{amsrefs}   
\usepackage[utf8]{inputenc}
\usepackage{tikz}
\usepackage{amsmath}
\usepackage{enumerate}
\usepackage{setspace}
\usepackage{lineno}
\usepackage{pgfplots}

\theoremstyle{plain}
\newtheorem{thm}{Theorem}[section]
\newtheorem{fact}[thm]{Fact}

\newtheorem{clm}[thm]{Claim}

\newtheorem{problem}[thm]{Problem}

\newtheorem{obs}[thm]{Observation}
\newtheorem*{clm*}{Claim}

\theoremstyle{definition}

\newtheorem{conj}[thm]{Conjecture}

\numberwithin{equation}{section}

\normalsize                              
\let\theta=\vartheta
\let\rho=\varrho
\let\phi=\varphi

\def\cP{{\mathcal P}}

\usepackage{enumitem} 

\DeclareMathOperator{\ex}{ex}

\let\polishlcross=\l
\def\l{\ifmmode\ell\else\polishlcross\fi}

\def\P{\mathcal{P}}

\theoremstyle{plain}
\let\epsilon=\varepsilon

\newtheoremstyle{note}% name
  {4pt}%      Space above
  {4pt}%      Space below 
  {\sl}%      Body font
  {}%         Indent amount (empty = no indent, \parindent = para indent)
  {\itshape}% Thm head font
  {.}%        Punctuation after thm head
  {.5em}%     Space after thm head: " " = normal interword space;
  {}%         Thm head spec (can be left empty, meaning `normal')

\theoremstyle{note}

\begin{document}

%\linenumbers 
%\onehalfspace

\title[]{Tur\'an number of disjoint triangles in 4-partite graphs}
\thanks{
Jie Han is partially supported by a Simons Collaboration Grant 630884.
Yi Zhao is partially supported by an NSF grant DMS 1700622 and Simons Collaboration Grant 710094.}
%
%%\author[W. Bedenknecht]{Wiebke Bedenknecht}
\author{Jie Han}
\address{Department of Mathematics, University of Rhode Island, 5 Lippitt Road, Kingston, RI, 02881, USA}
\email{jie\_han@uri.edu}
%%\author[Y. Kohayakawa]{Yoshiharu Kohayakawa}
\author{Yi Zhao}
%%\address{Fachbereich Mathematik, Universit\"at Hamburg,
%%  Bundesstra\ss{}e~55, D-20146 Hamburg, Germany}
%%\email{Wiebke.Bedenknecht@uni-hamburg.de}
%
%\address{Instituto de Matem\'atica e Estat\'{\i}stica, Universidade de
%	S\~ao Paulo, S\~ao Paulo, Brazil}
%\email{jhan@ime.usp.br}
%
\address
{Department of Mathematics and Statistics, Georgia State University, Atlanta, GA 30303}
\email{yzhao6@gsu.edu}

\maketitle

%\keywords{Hamiltonian cycle, random hypergraph, perturbed hypergraph}
%\subjclass[2010]{%05C35 (primary), 05C65, 05C80 (secondary)
%}
\begin{abstract}
Let $k\ge 2$ and $n_1\ge n_2\ge n_3\ge n_4$ be integers such that $n_4$ is sufficiently larger than $k$.
We determine the maximum number of edges of a 4-partite graph with parts of sizes $n_1,\dots, n_4$ that does not contain $k$ vertex-disjoint triangles.
For any $r> t\ge 3$,  we give a conjecture on the maximum number of edges of an $r$-partite graph that does not contain $k$ vertex-disjoint cliques $K_t$.
%We also determine the largest possible minimum degree among all $r$-partite triangle-free graphs.
\end{abstract}

%\tableofcontents

%We say a graph $G$ is $kK_r$-free if it does not contain $k$ vertex-disjoint copies of $K_r$.
%
%\begin{theorem}
%Given $k\ge 1$ and $r\ge 3$, let $n_1\ge n_2\ge \cdots \ge n_{r+1} \gg k$.
%Let $G$ be a $kK_r$-free $(r+1)$-partite graph with part sizes $n_1, \dots, n_{r+1}$.
%Then
%\begin{align*}
%\label{eq:ex}
%e(G)\le e(H) - \min\{&n_{r-2}n_{r+1}+n_{r-1}n_r - (k-1)n_{r-2}, \\
%&n_{r-1}n_r+n_{r-1}n_{r+1}+n_{r}n_{r+1}-(k-1)(n_{r-1}+n_{r})\},
%\end{align*}
%where $H$ is the complete $(r+1)$-partite graph with part sizes $n_1, \dots, n_{r+1}$.
%\end{theorem}
%
%\begin{proof}
%The case $r=3$ is Theorem~\ref{thm:main}.
%So we may assume that $r\ge 4$.
%We use the progressive induction on $n$ and the ordinary mathematical induction on $r$.
%We first assume that $n_1>n_{r-2}$.
%Then there exists $i\in [r-3]$, such that $n_i-1\ge n_{r-2}$.
%Since $n_i-1\ge n_{r-2}$, namely, for part sizes $n_1, \dots, n_{i-1}, n_i-1, n_{i+1} \dots, n_{r+1}$, the smallest four parts are still $n_{r-2},\dots, n_{r+1}$, then
%\[
%ex(kK_r; n_1, \dots, n_{r+1}) - ex(kK_r; n_1, \dots, n_{i-1}, n_i-1, n_{i+1} \dots, n_{r+1}) = \sum_{j\neq i} n_j.
%\]
%Therefore, for any vertex $v\in V_i$, if $d(v)<\sum_{j\neq i} n_j$, then we can delete $v$ and apply progressive induction.
%Otherwise we may assume that for any $v\in V_i$, $d(v)=\sum_{j\neq i} n_j$.
%Then we delete $V_i$ and apply induction on $r$, which is possible because every $v\in V_i$ connects to all vertices in other classes, and thus $G\setminus V_i$ must be $kK_{r-1}$-free.
%
%Now assume that $n_1 = n_{r-2}$.
%\end{proof}

\section{Introduction}  \label{sec:intro}

%Tur\'an-type problems are central in extremal graph theory.
%Mantel's theorem~\cite{mantel} is the very first result in extremal graph theory, which states that the maximum number of edges in an $n$-vertex triangle-free graph is $\lfloor n^2/4 \rfloor$, the number of edges of the balanced complete bipartite graphs.
%This was later generalized by Tur\'an~\cite{turan} to graphs without any $r$-cliques, $r\ge 3$.
Given two graphs $G$ and $F$, we say that $G$ is \emph{$F$-free} if $G$ does not contain $F$ as a subgraph. Let $K_t$ denote a complete graph on $t$ vertices, and $T_{n,t}$ denote a balanced complete $t$-partite graph on $n$ vertices (now known as the \emph{Tur\'an graph}).
In 1941, Tur\'an~\cite{turan} proved that $T_{n,t}$ has the maximum number of edges among all $K_{t+1}$-free graphs (the case $t=2$ was previously solved by Mantel~\cite{mantel}).
Tur\'an's result initiates the study of Extremal Graph Theory, an important area of research in modern Combinatorics (see the monograph of Bollob\'{a}s \cite{MR506522}).
%Tur\'an-type problems have a long and profound history.
%Note that Tur\'an problems become very hard in hypergraphs. For example, despite many efforts and recent developments, we still do not know the Tur\'an number for tetrahedron in a 3-uniform hypergraph.
Let $kK_t$ denote $k$ disjoint copies of $K_t$.
Simonovits~\cite{Simonovits} studied the Tur\'an problem for $kK_t$ and showed that when $n$ is sufficiently large, the (unique) extremal graph on $n$ vertices is the join of $K_{k-1}$ and the Tur\'an graph $T_{n-k+1, t-1}$.

In this paper we consider Tur\'an problems in multi-partite graphs. Let $K_{n_1,n_2,\dots, n_r}$ denote the complete $r$-partite graph on parts of sizes $n_1,n_2,\dots, n_r$.
This variant of the Tur\'an problem was first considered by Zarankiewicz~\cite{zaran}, who was interested in the case of forbidding $K_{s,t}$ in (subgraphs of) $K_{a,b}$.
Formally, given graphs $H$ and $F$, we define ex$(H, F)$ as the maximum number of edges in an $F$-free subgraph of $H$.
Bollob\'{a}s, Erd\H{o}s, and Straus \cite{MR379256} (see also \cite[Page 544]{MR506522}) proved the following result.
For any subset $I\subseteq [r]$, write $n_I := \sum_{i\in I} n_i$.

\begin{thm}\cite{MR379256}\label{thm:BES}
The extremal number $\ex(K_{n_1,\dots, n_r}, K_t)$ is equal to
\[
\max_{\P} \sum_{I\ne I'\in \mathcal P} n_I\cdot n_{I'},
\]
where the maximum is taken over all partitions $\mathcal P$ of $[r]$ into $t-1$ parts.
\end{thm}

The problem of forbidding disjoint copies of cliques in multi-partite graphs has been studied recently.
Chen, Li and Tu~\cite{CLT} determined $\ex(K_{n_1,n_2}, kK_2)$ %$\ex(K_{m,n}, kK_2)=m(k-1)$ for $1\le k\le n\le m$ 
and De Silva, Heysse and Young~\cite{DHY} showed that $\ex(K_{n_1, \dots, n_r}, kK_2) = (k-1)(n_1 + \dots + n_{r-1})$ for $n_1\ge \cdots \ge n_r$.
%extended it by replacing $K_{m, n}$ by an $r$-partite graph with $r\ge 3$.
De Silva, Heysse, Kapilow, Schenfisch and Young~\cite{DHKSY} determined $\ex(K_{n_1,\dots, n_r}, kK_r)$
and raised the question of determining $\ex(K_{n_1,\dots, n_r}, kK_t)$ when $r>t$.
After giving another proof of Theorem~\ref{thm:BES}, Bennett, English and Talanda-Fisher~\cite{MR3952133}  reiterated this question.
%namely, when the number of parts in the host graph is larger than the order of the forbidden clique.

\begin{problem}\cite{DHKSY}\label{pro1}
Determine $\ex(K_{n_1,\dots, n_r}, kK_t)$ when $r>t$.
\end{problem}

%The authors of~\cite{MR3952133} again asked for the solution of Problem~\ref{pro1} for $k\ge 2$.
In this paper we solve Problem~\ref{pro1} for $r=4$ and $t=3$ when all $n_i$'s are sufficiently large.
To state our result, for $k\ge 1$, we define a function of positive integers $n_1\ge n_2\ge n_3\ge n_4$: 
\begin{align*}
%\label{eq:ex}
g_k(n_1, n_2, n_3, n_4) &:=\max\left\{ (n_1+n_4)(n_2+n_3)+ (k-1)n_1, n_1(n_2+n_3+n_4) + (k-1)(n_2+n_3) \right\} \\
&= \left\{\begin{array}{lr}
         (n_1+n_4)(n_2+n_3) + (k-1) n_1 & \text{if } n_1\le n_2+n_3, \\
        n_1(n_2+n_3+n_4) + (k-1)(n_2+n_3),  & \text{if } n_1> n_2+n_3.
        \end{array}\right.
\end{align*}
When $G$ is a 4-partite graph with parts of sizes $n_1\ge n_2\ge n_3\ge n_4$, we define $g_k(G):= g_k(n_1, n_2, n_3, n_4)$. 
For arbitrary positive integers $a, b, c, d$, we define 
%$g_k(a, b, c, d)$ to be the function value corresponding to the non-ascending order of $a,b, c, d$. That is, 
$g_k(a,b,c,d)=g_k(a_1, a_2, a_3, a_4)$, where $a_1\ge a_2\ge a_3\ge a_4$ is a reordering of $a, b, c, d$.
%and $\{a_1, a_2, a_3, a_4\}=\{a,b,c,d\}$ as two multisets.

\begin{thm}\label{thm:main}
Given $k\ge 1$, there exists $N_0(k)$ such that if $G$ is a $kK_3$-free 4-partite graph with parts of sizes $n_1\ge n_2\ge n_3\ge n_4\ge 6k^2$ and $n_1+n_2+n_3+n_4\ge N_0(k)$, then $e(G)\le g_k(n_1, n_2, n_3, n_4) $.
In other words, $\ex(K_{n_1, n_2, n_3, n_4}, kK_3)\le g_k(n_1, n_2, n_3, n_4)$.
%, where
%\[
%%\label{eq:ex}
%g(n_1, n_2, n_3, n_4) :=\max\{ (n_1+n_4)(n_2+n_3)+ (k-1)n_1, n_1(n_2+n_3+n_4) + (k-1)(n_2+n_3) \}.
%\]
\end{thm}

Theorem~\ref{thm:main} is tight due to two constructions $G_1$ and $G_2$ below.
In fact, a subgraph of $G_2$ was given by De Silva et al. \cite{DHKSY} as a potential extremal construction; later Wagner~\cite{Wag19} realized that $G_1$ was a better construction for the $n_1 = n_2 = n_3 = n_4$ case.
%Let $G$ be a 4-partite with part sizes $n_1\ge n_2\ge n_3\ge n_4\ge k$.
Let $n_1\ge n_2\ge n_3\ge n_4\ge k$.
We define two 4-partite graphs with parts $V_1,\dots, V_4$ such that $|V_i|=n_i$.
Fix a set $Z$ of $k-1$ vertices in $V_4$.
Let
\[
G_1: =K_{V_1\cup V_4, \,V_2\cup V_3}\cup K_{Z, \,V_1} \text{ and } G_2:=K_{V_1, \,V_2\cup V_3\cup V_4}\cup K_{Z, \,V_2\cup V_3},
\] 
where $K_{V_1, \dots, V_r}$ denotes the complete $r$-partite graph with parts $V_1, \dots, V_r$. 
%be the union of two complete bipartite graphs, $K_{V_1, \,V_2\cup V_3\cup V_4}$ and $K_{Z, \,V_2\cup V_3}$.
Note that each triangle must intersect $Z$ and thus both $G_1$ and $G_2$ are $kK_3$-free.
Moreover, $e(G_1)=(n_1+n_4)(n_2+n_3)+ (k-1)n_1$ and $e(G_2)=n_1(n_2+n_3+n_4) + (k-1)(n_2+n_3)$.
%On the other hand, let $G_1$ be union of two complete bipartite graphs, $K_{V_1\cup V_4, \,V_2\cup V_3}$ and $K_{Z, \,V_1}$.
Thus $e(G_2)\le e(G_1)$ if and only if $n_1\le n_2+n_3$ and equality holds when $n_1=n_2+n_3$.

\begin{figure}[h]
\begin{center}
\begin{tikzpicture}
[inner sep=2pt,
%   vertex/.style={circle, draw=blue!50, fill=blue!50},
   vertex3/.style={circle, inner sep=7,minimum size=1.7cm},
   vertex2/.style={circle, inner sep=7,minimum size=1.8cm},
   vertex1/.style={circle, inner sep=7,minimum size=1.9cm},
   vertex/.style={circle, inner sep=7,minimum size=2cm},
   ver/.style={circle, minimum size=0.01, fill=blue!20},
   ]
%\begin{pgfonlayer}{bg}    % select the background layer
%\draw[rounded corners] (-1,0) rectangle (9, 1.8);
%\draw[rounded corners] (-1,-2) rectangle (9, -0.2);
%\end{pgfonlayer}

%
\node[vertex, draw] (v1) at (0,2.15) {$V_1$};
\node[vertex2, draw] (v3) at (3,0) {$V_3$};
\node[vertex1, draw] (v2) at (0,0) {$V_2$};
\node[vertex3, draw] (v4) at (3,2) {$V_4$};

\node[vertex, draw] (V1) at (7,2.15) {$V_1$};
\node[vertex2, draw] (V3) at (10,0) {$V_3$};
\node[vertex1, draw] (V2) at (7,0) {$V_2$};
\node[vertex3, draw] (V4) at (10,2) {$V_4$};

%\node (LL) at (4.6, -0.1) {$L'$};
\draw[line width=0.2cm]  (v1) -- (v2); 
\draw[line width=0.2cm]  (v1) -- (v3); % S-S_1
\draw[line width=0.2cm]  (v4) -- (v2); % S-S_1
\draw[line width=0.2cm]  (v4) -- (v3); % S-S_1

\draw[line width=0.2cm]  (V1) -- (V2); % S-S_1
\draw[line width=0.2cm]  (V1) -- (V3); % S-S_1
\draw[line width=0.2cm]  (V1) -- (V4); % S-S_1

\node[ver, draw] (z) at (2.5,2) {$Z$};
\draw[line width=0.15cm, color=blue!50]  (v1) -- (z); 
\node[ver, draw] (z2) at (9.5,1.7) {$Z$};
\draw[line width=0.15cm, color=blue!50]  (V2) -- (z2); 
\draw[line width=0.15cm, color=blue!50]  (V3) -- (z2); 

\end{tikzpicture}

\caption{The extremal graphs $G_1$ and $G_2$}
\end{center}
\end{figure}

Our proof uses a \emph{progressive induction} (an induction without a base case) on the total number of vertices and a standard induction on~$k$ that uses Theorem~\ref{thm:BES} as the base case.

\medskip
We conjecture an answer to Problem~\ref{pro1} in general, which includes all aforementioned results \cite{MR3952133, CLT, DHY} and Theorem~\ref{thm:main}.
\begin{conj}
\label{conj}
Given $r> t\ge 3$ and $k\ge 2$, let $n_1, \dots, n_r$ be sufficiently large.
For $I\subseteq [r]$, write $m_I := \min_{i\in I}n_i$.
Given a partition $\mathcal P$ of $[r]$, let $n_{\mathcal P}:=\max_{I\in \mathcal P}\{n_I-m_I\}$.
The Tur\'an number $\ex(K_{n_1,\dots, n_r}, kK_t)$ is equal to
\begin{equation}\label{eq:conj}
\max_{\P} \left\{ (k-1)n_{\mathcal P} + \sum_{I\ne I'\in \mathcal P} n_I\cdot n_{I'} \right\},
\end{equation}
where the maximum is taken over all partitions $\mathcal P$ of $[r]$ into $t-1$ parts.
\end{conj}

The bound~\eqref{eq:conj} is achieved by the following graph.
Given integers $k, t$ and $n_1,\dots, n_r$ with $r>t$ and $n_i\ge k$ for all $i$, let $\cP$ be a partition of $[r]$ into $t-1$ parts that maximizes~\eqref{eq:conj}.
Let $G$ be an $r$-partite graph whose parts have sizes $n_1,\dots, n_r$.
Partition $G$ into $t-1$ parts according to $\cP$, namely, let $V_{I}=\bigcup_{i\in I}V_i$ for every $I\in \cP$ and include all edges between $V_I$ and $V_{I'}$ for all $I\ne I'\in \cP$.
In addition, let $I_0\in \mathcal P$ such that $n_{\cP}=n_{I_0}-m_{I_0}$ and let $V_{i_0}$ be the smallest part in $V_{I_0}$.
We choose a set $Z\subseteq V_{i_0}$ of $k-1$ vertices and add all edges between $Z$ and $V_{I_0}\setminus V_{i_0}$.

Verifying Conjecture~\ref{conj} seems hard due to the complexity of \eqref{eq:conj} -- we shall discuss this in the last section.

\medskip
\noindent\textbf{Notation.}
Given a graph $G=(V, E)$, let $|G|$ denote the order of $G$. Suppose $A, B$ are two disjoint subsets of $V$.
%, let $e(G)$ denote its number of edges.
%For disjoint vertex sets $A, B\subseteq V$, l
Let $e(A):=e(G[A])$ be the number of edges of $G$ in $A$ and $e(A, B)$ be the number of edges of $G$ with one end in $A$ and the other in $B$.
Moreover, let $G\setminus A:=G[V\setminus A]$.
Denote by 
\[
e(A; G) := e(G) - e(G\setminus A), 
\]
the number of edges of $G$ incident to $A$. 
Given a vertex $x$, let $N(x)$ denote the set of neighbors of $x$.
For vertices $x, y$ and $z$, we often write $xyz$ for $\{x, y, z\}$.
We sometimes abuse this notation by using $xy\in A\times B$ to indicate that $x\in A$ and $y\in B$. %(then $xy$ represents an ordered pair here).
Given an $r$-partite graph $G$, a \emph{crossing set} is a set that contains at most one vertex from each part of $G$.

\section{Proof of Theorem~\ref{thm:main}}

In this section we prove Theorem~\ref{thm:main}. 
Define two sequences $N_0(k)$ and $M_0(k)$ recursively by letting $N_0(1)=1$, 
\begin{equation}\label{eq:M0}
M_0(k)=\max\{72(k-1)^3, 96k^2, N_0(k-1)+3\}, \quad \text{and} \quad  N_0(k)= M_0(k)^2
\end{equation}
for $k\ge 2$.
Given a 4-partite graph $G$, let  $v_4(G)$ denote the size of the smallest part of $G$.
Define $\phi(G) := e(G) - g_k(G)$.
The following theorem is the main step in the proof of Theorem~\ref{thm:main}. 
\begin{thm}
\label{thm:key}
Suppose $k\ge 2$ and Theorem~\ref{thm:main} holds for $k-1$. 
Let $G$ be a 4-partite graph of order $|G| > M_0(k)$ and with $v_4(G) \ge 6k^2$. If $G$ is $k K_3$-free and $\phi(G)> 0$, then we can find a subgraph $G'$ of $G$ such that $|G| -2 \le |G'| \le |G| - 1$, $v_4(G')\ge 6k^2$, and $\phi(G') > \phi(G)$. 
\end{thm}

Theorem~\ref{thm:main} nows follows from Theorem~\ref{thm:key} by an induction on $k$ and a progressive induction on $|G|$ (e.g., used in \cite{Simonovits}). 

\begin{proof}[Proof of Theorem~\ref{thm:main}]
The base case $k=1$ follows from Theorem~\ref{thm:BES} with $N_0(1)=1$. Let $k\ge 2$ and $G$ be a 4-partite graph of order $|G| \ge N_0(k)$ and with $v_4(G) \ge 6k^2$. Suppose $G$ is $k K_3$-free and $\phi(G)> 0$. By Theorem~\ref{thm:key}, we find a subgraph $G_1\subset G$ such that $|G| -2 \le |G_1| \le |G| - 1$, $v_4(G_1)\ge 6k^2$, and $\phi(G_1) > \phi(G)\ge 1$. Repeating this process, we obtain subgraphs $G_1\supset G_2\supset G_3 \supset \cdots \supset G_t$ such that $|G| - 2i \le |G_i|\le |G| - i$ and $\phi(G_i) > i$ for $i=1, \dots, t$. We stop at $G_t$ because $|G_t|\le M_0(k)$. Hence,
\[
t \ge \frac{|G|- |G_t|}2 \ge  \frac{N_0(k)- M_0(k)}2 = \frac{M_0(k)^2 - M_0(k)}2 = \binom{M_0(k)}{2}.
\]
Consequently, $\phi(G_t) > \binom{M_0(k)}{2}$. However, this is impossible because $\phi(G_t)\le e(G_t) \le  \binom{M_0(k)}{2}$.
\end{proof}

%We thus require $N_0(k)\ge M_0^2(k)$, where $M_0(k)$ is defined in \eqref{eq:M0}; in particular, $N_0(k)\ge (N_0(k-1)+3)^2$. It is easy to see that $N_0(k)=2^{3^k}$ suffices.

The rest of this section is devoted to the proof of Theorem~\ref{thm:key}. 

\begin{proof}[Proof of Theorem~\ref{thm:key}]
Let $k\ge 2$ and suppose that 
\begin{enumerate}[label=($\ast$)]
\item for any $(k-1)K_3$-free 4-partite graph $\tilde{G}$ with part sizes $n_1'\ge n_2'\ge n_3'\ge n_4'\ge 6(k-1)^2$ and $\sum_{i\in [4]}n_i' \ge N_0(k-1)$, we have $e(\tilde{G})\le g_{k-1}(n_1', n_2', n_3', n_4')$.
\label{item:IH}
\end{enumerate}

%To prove the assertion for $k$, we apply  Lemma~\ref{lem:prog}. For $N\ge 24k^2$, we let $\fA_N$ be the collection of all $kK_3$-free 4-partite graphs whose parts have sizes $n_1\ge \cdots \ge n_4$ such that $N=n_1+\cdots+n_4$ and $n_4\ge 6k^2$; let $\fA_N= \emptyset$ for $N<24k^2$.
%% for some function $f$ of $k$.  
%Let $\mathfrak{A}=\bigcup_{N\in \mathbb N} \fA_{N}$.
%Define $B$ as the property that $e(G)\le g_k(n_1, n_2, n_3, n_4)$ and $\phi(G) := \max\{0, e(G) - g_k(n_1, n_2, n_3, n_4)\}$.
%Thus Lemma~\ref{lem:prog}~\ref{item:a} holds and it remains to prove~\ref{item:b}. 
%
%Let $N>M_0(k)$, where
%\begin{equation}\label{eq:M0}
%M_0(k):=\max\{432 k^4/(k-1)+ 6k^2, N_0(k-1)+3\}.
%\end{equation}
%%Note that $M_0(k) - 3 = 2^{3^{k-1}} = N_0(k-1)$.

Let $G$ be a 4-partite graph of order $|G| > M_0(k)$ and with parts of size $n_1\ge n_2\ge n_3\ge n_4 \ge 6k^2$. Assume that $G$ is $kK_3$-free and $\phi(G)> 0$. Without loss of generality, we assume that $G$ contains $k-1$ disjoint triangles -- otherwise we keep adding edges to $G$ until it contains $k-1$ disjoint triangles (as a result, $\phi(G)$ increases).
%By adding edges to $G$, we may also assume that $G$ is the densest $kK_3$-free 4-partite graph on parts $$n_1, n_2, n_3, n_4$ and $\phi(G)> 0$. Hence 
%Let $G_0$ does not satisfy $B$ (thus $\phi(G_0)>0$).
%Let $G$ be the largest $kK_3$-free 4-partite supergraph of $G_0$ on the same parts of $G_0$. Since adding any additional edge to $G$ would create a copy of $kK_3$, it implies that $G$ contains $k-1$ disjoint triangles.
%Clearly $0<\phi(G_0)\le \phi(G)$.
Our goal is to show that there exists a crossing set $T\subset V(G)$ of size at most 2 such that $\phi(G) < \phi(G\setminus T)$ and $v_4(G\setminus T)\ge 6k^2$.
%\begin{enumerate}[label=($\dagger$)]
%\item there exists some crossing set $T\subset V(G)$ of size at most 2 such that $\phi(G) < \phi(G\setminus T)$ and all parts of $G\setminus T$ have sizes at least $6k^2$, \label{item:bbb}
%\item $\phi(G)=0$,
%\end{enumerate}
%which would prove~\ref{item:b} because $G\setminus T\in \fA_{N'}$, $N/2 < N-2\le N'\le N-1$ and $\phi(G_0)\le \phi(G) < \phi(G\setminus T)$.

%In the actual proof, we also use an induction on $k$, where the base case $k=1$ follows from Theorem~\ref{thm:BES}.
%Note that the size requirement will always be satisfied, as we only apply the inductive hypothesis on an induced subgraph of $G$ with at most three vertices removed.
%This is valid because $n_1+n_2+n_3+n_4\ge M_0(k) \ge N_0(k-1)+3$, and $n_4\ge 6k^2 \ge 6(k-1)^2+3$.

We proceed in the following cases. It is easy to see that these cases cover all possibilities. In each case we verify $v_4(G\setminus T)\ge 6k^2$ immediately.

%Since $n_4\ge 6k^2$, in each case it suffices to show that $\phi(G) < \phi(G\setminus T)$ and

\noindent\textbf{Case 0. $n_1 > n_2 + n_3$.} We will select a one-element set $T\subset V_1$. Since $n_1> 2 n_4$, we have $n_1 - 1> n_4$ and thus $v_4(G\setminus T)= n_4 \ge 6k^2$.

We assume $n_1\le n_2+n_3$ in the remaining cases.%\footnote{YZ: it seems only Case 1 may use $|T|=2$ and Case 4 only uses $T\subset V_1$.}

\noindent\textbf{Case 1. $n_1 > n_3$ and $n_2> n_4$.} We will select a crossing set $T\subset V_1\cup V_2$. Since $n_1 - 1 \ge n_2 - 1 \ge n_4$, we have $v_4(G\setminus T) = n_4  \ge 6k^2$.

\noindent\textbf{Case 2. $n_1 = n_2 = n_3 \ge n_4 > 6k^2$.} We select a one-element set $T\subset V(G)$. Then $v_4(G\setminus T)\ge n_4-1\ge 6k^2$.

\noindent\textbf{Case 3. $n_1 = n_2 = n_3 > n_4 = 6k^2$.} We will select a one-element set $T\subset V_1\cup V_2\cup V_3$. Since $n_3 -1 \ge n_4$, we have $v_4(G\setminus T) = n_4  = 6k^2$.

\noindent\textbf{Case 4. $n_1 > n_2 = n_3 = n_4 $.} We will select a one-element set $T\subset V_1$. Since $n_1>n_4$, $v_4(G\setminus T)= n_4\ge 6k^2$.

%\noindent\textbf{Case 5. $n_1 \ge n_2 = n_3 = n_4 = 6k^2$.} Since $n_1\le n_2+n_3$, we conclude that $n_1+n_2+n_3+n_4\le 30k^2 < M_0$ so this case does not exist.

%Now suppose under each case we can show our promise, then the proof is completed because $n_4\ge 6k^2$ is maintained throughout the proof (crucially, in \textbf{Cases 0, 1} and \textbf{3}, $T\cap V_4=\emptyset$, so the resulting $n_4$ will not decrease).
%Note that this finishes the proof as $n_4>6k^2$ except in \textbf{Cases 0, 1} and \textbf{3}, in which cases $T\cap V_4=\emptyset$ so that the smallest part size of $G\setminus T$ is still at least $6k^2$.
%The bound $N=2^{3^k}$ can be shown by going through the proof of Lemma~\ref{lem:prog} and noting that we actually remove at most two vertices in each round.
%So when the condition~\ref{item:b} stops working or when we arrive Case 5, the resulting graph $G^*$ has at most $M_0=2^{3^{k-1}}+3$ vertices, namely, the progressive induction has proceeded at least
%\[
%(N-M_0)/2 > 2^{3^k-1} > (2^{3^{k-1}}+3)^2 = M_0^2
%\]
%rounds, and thus the resulting graph satisfies $\phi(G^*)\ge (N-M_0)/2 > M_0^2$, which is impossible.
%The proof is completed (modulo the validation of the cases).

\medskip
It remains to show $\phi(G) < \phi(G\setminus T)$ in \textbf{Cases 0--4}. This is actually easy in \textbf{Case 0}.
%In fact, \textbf{Case 0} can be proved easily.

%\begin{proof}[Proof of {\rm \textbf{Case 0}}]
\noindent \textbf{Case 0.}
Recall that $\phi(G) = e(G) - g_k(n_1, n_2, n_3, n_4) > 0$. Since $n_1 > n_2 + n_3$, 
\[
g_k(n_1, n_2, n_3, n_4) = n_1(n_2+n_3+n_4) + (k-1)(n_2+n_3).
\]
First assume that some vertex $v\in V_1$ satisfies $d(v)<n_2+n_3+n_4$. Let $T=\{v\}$. Since $n_1-1\ge n_2+n_3$, %we have
\begin{align*}
g_k(n_1-1, n_2, n_3, n_4) &=  (n_1 - 1)(n_2+n_3+n_4) + (k-1)(n_2+n_3)\\
&= g_k(n_1, n_2, n_3, n_4) - (n_2+n_3+n_4).
\end{align*}
It follows that 
\[
\phi(G\setminus \{v\}) = e(G) - d(v) - g_k(n_1-1, n_2, n_3, n_4) > e(G) - g_k(n_1, n_2, n_3, n_4) = \phi(G),
\]
%\begin{align*}
%g_k(n_1, n_2, n_3, n_4) - g_k(n_1-1, n_2, n_3, n_4) &= n_1(n_2+n_3+n_4) + (k-1)(n_2+n_3) - \\
%& \left( (n_1 - 1)(n_2+n_3+n_4) + (k-1)(n_2+n_3) \right)\\
%&= n_2+n_3+n_4.
%\end{align*}
as desired. Otherwise, $G[V_1, V_2\cup V_3\cup V_4]$ must be complete. Since $G$ is $k K_3$-free, it follows that $G[V_2\cup V_3\cup V_4]$ contains no matching of size $k$. The result of ~\cite{DHY} or a simple induction on $k$\footnote{If there is a vertex of degree at least $2k-1$, then we can delete it and apply induction; otherwise, as the size of the maximum matching is $k-1$, there are at most $2(k-1)(2k-1)\le (k-1)(n_2 + n_3)$ edges (using $k\ll n_3\le n_2$).}  yields that $e(G[V_2\cup V_3\cup V_4]) \le (k-1)(n_2+n_3)$.
This shows that $e(G)\le n_1(n_2+n_3+n_4) + (k-1)(n_2+n_3)$, namely, $\phi(G)=0$, a contradiction. 
%\end{proof}

\medskip
In the rest of the proof we assume $n_1\le n_2+n_3$ and will resolve \textbf{Cases 1--4}.
%In order to derive \ref{item:b}, we consider a set $T\subseteq V(G)$ of size at most two and let $n'_1, n'_2, n'_3, n'_4$ denote the sizes of the parts of $G\setminus T$.
%%Then $G\setminus T\in \fA_{n'_1+\cdots+ n'_4}$.
%Then $\phi(G) < \phi(G\setminus T)$ is equivalent to 
%\[
%e(G) -  g_k(n_1, n_2, n_3, n_4) < e(G\setminus T) - g_k(n'_1, n'_2, n'_3, n'_4).
%\]
%%Denote by $e(T; G) := e(G) - e(G\setminus T)$, the number of edges of $G$ incident to $T$. 
%Then $\phi(G) < \phi(G\setminus T)$ is equivalent to $e(T; G)< g_k(n_1, n_2, n_3, n_4) - g_k(n'_1, n'_2, n'_3, n'_4)$.
%Thus, when proving Theorem~\ref{thm:main} by contradiction, we may assume that for all sets $T\subseteq V(G)$ of size at most two,
%\begin{equation}\label{eq:etg}
%e(T; G)\ge g_k(n_1, n_2, n_3, n_4) - g_k(n'_1, n'_2, n'_3, n'_4).
%\end{equation}
%

%\section{4-partite $kK_3$-free graphs}
One difficulty in these cases is that, after we delete a set $T\subseteq V(G)$, the sizes of the four parts of $G\setminus T$ may not follow the order in $G$.
For instance, suppose $n_1\le n_2 + n_3$ and $T=\{v\}\subseteq V_1$.
If $n_1>n_2$, then  the order of the part sizes of $G\setminus T$ is $n_1-1\ge n_2\ge n_3\ge n_4$, the same as in $G$.
%By~\eqref{eq:etg}, we may assume that for every $v\in V_1$,
%\[
%d(v)=e(T; G)\ge g_k(n_1, n_2, n_3, n_4) - g_k(n_1-1, n_2, n_3, n_4) =n_2+n_3+k-1,
%\]
%which matches the degree of the vertices in $V_1$ in the extremal graph $G_1$.
%and the $g$ function well.
%(Indeed, if $n_1>n_2> n_3>n_4$, then deleting a set $T$ with at most one vertex from each class will not affect the order of the parts.)
However, when $n_1=n_2>n_3\ge n_4$, the order of the part sizes of $G\setminus T$ is $n_2\ge n_1-1\ge n_3\ge n_4$, and the degree estimates we obtain are quite different.
%\footnote{YZ: why? I don't think so. In fact, if we are allowed relabeling $V_i$, then the inequality should be the same.}
%and we can only derive $d(v)\ge n_1+n_4$ from~\eqref{eq:etg}. 
Another complication comes from the fact that there are two possible extremal graphs. Even under the assumption that $n_1 \le n_2 + n_3$, we still have to consider the possibility of $n'_1> n'_2 + n'_3$ in $G\setminus T$, where $n'_1, n'_2, n'_3, n'_4$ are the part sizes of $G\setminus T$.

Although a case analysis is inevitable, we study the structure of $G$ in Section ~\ref{sec:21} and use it to simplify the presentation of the proofs of \textbf{Cases 1--4} in Section~\ref{sec:22}.

\subsection{Preparation}
\label{sec:21}
%Theorem~\ref{thm:main} is proved in this section.
%Suppose $N$ is sufficiently large and in particular
%\begin{equation}\label{eq:n4}
%n_4\ge N\ge 4k.
%\end{equation}
%We shall use a progressive induction on $n_1+\cdots+ n_4$ and a standard induction on $k$.
%Let $k\ge 2$ as the case $k=1$ has been proved in Theorem~\ref{thm:BES}.
%We also assume that $G$ is maximal, that is, if we add any additional edge to $G$, then $kK_3\subseteq G$.
%%By this we can assume that $(k-1)K_3\subseteq G$.
%So $G$ contains at least $k-1$ disjoint triangles.

%Our main tool is the progressive induction on $n_1+\cdots+ n_4$ and induction on $k$.
%The case $n_1\le n_2+n_3$ is more involved and w
We first give several preliminary results.
%Consider a set $T\subseteq V(G)$ of size at most two and let $n'_1, n'_2, n'_3, n'_4$ denote the sizes of the parts of $G\setminus T$.
%%Then $G\setminus T\in \fA_{n'_1+\cdots+ n'_4}$.
%Then $\phi(G) < \phi(G\setminus T)$ is equivalent to 
%\[
%e(G) -  g_k(n_1, n_2, n_3, n_4) < e(G\setminus T) - g_k(n'_1, n'_2, n'_3, n'_4),
%\]
%which in turn is equivalent to $e(T; G)< g_k(n_1, n_2, n_3, n_4) - g_k(n'_1, n'_2, n'_3, n'_4)$.
%Thus, in each case, we may assume that $\phi(G) \ge \phi(G\setminus T)$ (as otherwise we are done), that is, for the sets $T\subseteq V(G)$ of size at most two as specified in each of \textbf{Cases 1 - 4},
%\begin{equation}\label{eq:etg}
%e(T; G) \ge (n_1+n_4)(n_2+n_3) + (k-1)n_1 - g_k(n'_1, n'_2, n'_3, n'_4) .
%\end{equation}
%
%Assume that $e(G)> g_k(n_1, n_2, n_3, n_4)$.
An edge of $G$ is called \emph{rich} if it is contained in at least $k$ triangles whose third vertices are located in the same part of $V(G)$.
We show that every triangle in $G$ must contain a rich edge and $G$ contains at most $6(k-1)^2$ rich edges.
Let $Z$ be the set of vertices incident to at least one rich edge.
Thus, not only is $G\setminus Z$ triangle-free,
%\footnote{YZ removed (which would be true for any $Z$ that contains a copy of $(k-1)K_3$). Does the referee want this sentence? It is almost impossible for $Z$ to contain a copy of $(k-1)K_3$. JH. I put it to emphasize the function of $Z$. Yeah I don't mind if it is removed.}
but also \emph{every edge in $G\setminus Z$ is not contained in any triangle of $G$} because such a triangle would not contain any rich edge.
%\footnote{YZ removed ``Then we derive a contradiction by counting the edges of $G$." because it does not fit here.}

We shall use the following simple fact.
\begin{fact}\label{fact:com_neigh}
Let $G$ be a 4-partite graph with parts $V_1,\dots, V_4$ and suppose $x\in V_1$ and $y\in V_2$.
Let $n_i:=|V_i|$ for $i\in [4]$.
Then $x$ and $y$ have at least $d(x)+d(y) - \sum_{i\in [4]}n_i$ common neighbors in $G$.
%In particular, if $x$ and $y$ have no common neighbor, then $d(x)+d(y) \le \sum_{i\in [4]}n_i$ and the equality holds if and only if $xy\in E(G)$, $V_2\subseteq N(x)$ and $V_1\subseteq N(y)$.
In particular, if $x$ and $y$ have no common neighbor, then $d(x)+d(y) = \sum_{i\in [4]}n_i$ implies that $xy\in E(G)$, $V_2\subseteq N(x)$ and $V_1\subseteq N(y)$.
Moreover, if $d(x)+d(y)\ge \sum_{i\in [4]}n_i+2k-1$, then $xy$ is rich.
\end{fact}

\begin{proof}
Note that $|N(x)\cap (V_3\cup V_4)| = d(x) - |N(x)\cap V_2| \ge d(x) - n_2$ and $|N(y)\cap ( V_3\cup V_4)|=d(y) - |N(y)\cap V_1|\ge d(y) - n_1$.
Let $m$ denote the number of common neighbors of $x$ and $y$.
Then $m\ge |N(x)\cap (V_3\cup V_4)| + |N(y)\cap (V_3\cup V_4)| - n_3 - n_4 \ge d(x)+d(y) - \sum_{i\in [4]}n_i$.
So the first part of the fact follows.
In particular, if $m=0$, then $d(x)+d(y) \le \sum_{i\in [4]}n_i$.
Moreover, if the equality holds, then the inequalities in previous calculations must be equalities.
In particular,  $V_2 \subseteq N(x)$ and $V_1\subseteq N(y)$, which also imply that $xy\in E(G)$.

For the ``moreover'' part, note that $d(x)+d(y)\ge \sum_{i\in [4]}n_i+2k-1$ implies that $x$ and $y$ have at least $2k-1$ common neighbors and thus at least  $k$ common neighbors in one part. Therefore $xy$ is rich.
\end{proof}

Recall that we have assumed that $\phi(G)>0$ and $n_1\le n_2 + n_3$. Thus, 
\begin{equation}
\label{eq:contra}
e(G)> g_k(n_1, n_2, n_3, n_4) = (n_1+n_4)(n_2+n_3)+ (k-1)n_1.
\end{equation}

Let $R$ be the subgraph of $G$ induced by the rich edges of $G$, and let $Z= V(R)$ be the set of the vertices of $G$ that are incident to at least one rich edge.
%We have the following claim. 
\begin{clm}
Suppose \ref{item:IH}, \eqref{eq:contra}, and $G$ is $k K_3$-free. Then the following assertions hold:
\begin{enumerate}[label=$(\roman*)$]
\item every vertex is contained in at most $k-1$ edges of $R$ whose other ends are located in the same part of $G$; in particular, the maximum degree of $R$ is at most $3k-3$; \label{item:1}
%\item $R$ does not have a matching of size $k$;
\item $e(R)\le 6(k-1)^2$ and $|Z|\le 6(k-1)^2$; \label{item:2}
\item every triangle in $G$ contains an edge in $R$. \label{item:3}
%\item for every edge $xy\in V_i\times V_j$ in $R$, $\{i,j\}=\{1,4\}$ or $\{2,3\}$, we have $d(x)+d(y) \ge 2n_t + 2n_l + n_4/2$, where $\{i,j,t,l\}=[4]$. \label{item:4}
\end{enumerate}
\end{clm}

\begin{proof}
We first show \ref{item:1} $\Rightarrow$ \ref{item:2}. Note that if $R$ has a matching of size $k$, then we can greedily build $k$ vertex-disjoint triangles by extending each rich edge in the matching. This contradicts the assumption that $G$ is $k K_3$-free. Therefore, the largest matching in $R$ is of size at most $k-1$ and consequently, $R$ has a vertex cover of size at most $2(k-1)$. If the maximum degree of $R$ is at most $3k-3$, then 
$e(R)\le 2(k-1)(3k-4) + k-1 < 6(k-1)^2$ and $|Z|\le 2(k-1)(3k-4) + 2(k-1) = 6(k-1)^2$, confirming \ref{item:2}.

To see \ref{item:1}, we assume that some vertex $v$ is incident to $k$ rich edges whose other ends are in the same part of $G$.
%then we will delete this vertex and proceed by the induction on $k$.
If there is a copy $S$ of $(k-1)K_3$ in $G\setminus \{v\}$, then we can pick a rich edge in $G\setminus S$ that contains $v$ and then extend this rich edge to a triangle that does not intersect $S$.
This gives a $kK_3$ in $G$, a contradiction.
Thus, we infer that $G\setminus \{v\}$ is $(k-1)K_3$-free.
%\footnote{YZ: Below needs an explanation. JH. I tried but now sure whether this is good or not. YZ gave a detailed explanation.} 

Let $n'_1 \ge n'_2 \ge n'_3 \ge n'_4$ be the sizes of four parts of $G\setminus \{v\}$. 
By \ref{item:IH}, we have $e(G\setminus \{v\})\le g_{k-1}(n'_1, n'_2, n'_3, n'_4)$.
To estimate $g_{k-1}(n'_1, n'_2, n'_3, n'_4)$, we first observe that there exists $i_0\in [4]$ such that $n'_i = n_i$ for all $i\ne i_0$ and $n_{i_0} = n_{i_0} - 1$; and furthermore, $n'_i = |V_i\setminus \{v\}|$ for $i\in [4]$ after relabeling $V_1, V_2, V_3, V_4$ if necessary (but maintaining $n_i = |V_i|$). This is obvious when $v\in V_{i_0}$ and $n_{i_0} > n_{i_0 + 1}$. Otherwise, for example, assume that $v\in V_1$ and $n_1 = n_2> n_3$ (other cases are similar). Then $n'_1 = n_2 = n_1$ and $n'_2 = n_1 -1 = n_2 -1$. After relabeling $V_1$ and $V_2$, we have $v\in V_2$, and $n'_i = |V_i\setminus \{v\}|$ for $i\in [4]$.

By the definition of $g$, we consider two cases. When $n'_1 \le n'_2 + n'_3$, we have
\begin{align}
g_{k-1}(n'_1, n'_2, n'_3, n'_4) &=  (n'_1+n'_4)(n'_2+n'_3)+ (k-2)n'_1 \nonumber \\
&\le  \left\{\begin{array}{lr}
         (n_1+n_4-1)(n_2+n_3) + (k-2)n_1 & \text{if } v\in V_1\cup V_4, \\
        (n_1+n_4)(n_2+n_3-1) + (k-2)n_1,  & \text{if } v\in V_2\cup V_3.
        \end{array}\right.
        \label{eq:23i}
\end{align}
Together with \eqref{eq:contra} and \ref{item:IH}, this implies that 
\begin{align*}
d_G(v) &= e(G) - e(G\setminus \{v\}) > g_k(n_1, n_2, n_3, n_4) - g_{k-1}(n'_1, n'_2, n'_3, n'_4) \\\
&\ge
  \left\{\begin{array}{lr}
        n_1+ n_2+n_3  & \text{if } v\in V_1\cup V_4, \\
        2n_1+n_4,  & \text{if } v\in V_2\cup V_3,
        \end{array}\right.
\end{align*}
which is impossible.
When $n'_1 > n'_2 + n'_3$, it must be the case when $n_1 = n_2 + n_3$ and $n'_{i_0} = n_{i_0} - 1$ for $i_0\in \{2, 3\}$. Thus 
\begin{align*}
g_{k-1}(n'_1, n'_2, n'_3, n'_4) &=  n'_1(n'_2+n'_3+ n'_4)+ (k-2)(n'_2 + n'_3)\\
&=  (n_2+n_3)(n_1+n_4-1) + (k-2)(n_1-1). 
\end{align*}
Together with \eqref{eq:contra} and \ref{item:IH}, this implies that $d_G(v)> n_1 + n_2 + n_3$, which is impossible for any $v\in V(G)$.
%Thus \eqref{eq:23i} holds for both cases.
%So we can bound $e(G\setminus \{v\})$ by the inductive hypothesis.
%There are four cases depending on the order of the part sizes of $G\setminus \{v\}$, but we always have the following upper bound Therefore
%\begin{align*}
%e(G\setminus \{v\}) \le \max\{ (n_1&+n_4-1)(n_2+n_3)+ (k-2)n_1, \\
%&(n_1+n_4)(n_2+n_3-1)+ (k-2)n_1 \},
%\end{align*}
%where the two cases depend on the location of $v$ and we give up on the possible saving of an additional $k-2$.
%where in the maximum, the former term is for the case $v\in V_1$
%Moreover, for the maximum above, if $v\in V_4$, then the first term must achieve the maximum, and $d_G(v)\le n_1+n_2+n_3$; otherwise, we have $d_G(v)\le n_1+n_2+n_4$.
%It is easy to check that $e(G) = e(G\setminus \{v\}) + d(v) \le (n_1+n_4)(n_2+n_3)+ (k-1)n_1$ holds for all cases, contradicting~\eqref{eq:contra}.
%Since $d_G(v) = e(G) - e(G\setminus \{v\}) > g_k(n_1, n_2, n_3, n_4) - g_{k-1}(n'_1, n'_2, n'_3, n'_4)$

To see \ref{item:3}, let $S$ be a triangle in $G$ and consider $G\setminus S$.
%\footnote{YZ: I suggest to preserve $T$ for those in the four cases and use $S$ here. Also, $e(T;G)$ has not been defined.}
Since $G$ is $k K_3$-free, $G\setminus S$ is $(k-1) K_3$-free.
By~\ref{item:IH}, we have $e(G\setminus S)\le g_{k-1}(n_1', n_2', n_3', n_4')$ where $n_1'\ge n_2'\ge n_3'\ge n_4'$ are the sizes of parts of $G\setminus S$.
We observe that there exists $i_0\in [4]$ such that $n'_i = n_i - 1$ for $i\ne i_0$ and $n'_{i_0} = n_{i_0}$; furthermore, $n'_i = |V_i\setminus S|$ after relabeling $V_1, V_2, V_3, V_4$ if necessary (while maintaining $n_i = |V_i|$). This is obvious when $S\subset \bigcup_{i\ne i_0} V_i$ and either $i_0=1$ or $n_{i_0 - 1}> n_{i_0}$. Otherwise, for example, assume that $S\subset V_1\cup V_2\cup V_3$ and $n_2 > n_3 = n_4$ (other cases are similar).  We have $n'_1 = n_1 -1$, $n'_2 = n_2 - 1$, $n'_3 = n_4 = n_3$ and $n'_4 = n_3 -1 = n_4 -1$. After swapping $V_3$ and $V_4$, we have $S\subset V_1\cup V_2\cup V_4$.

If $n_1' \le n_2'+n_3'$. then $g_{k-1}(n_1', n_2', n_3', n_4') = (n_1'+n_4')(n_2'+n_3')+ (k-2)n_1'$. By our observation on the values of $n'_1, n'_2, n'_3, n'_4$, it follows that 
\[
g_{k-1}(n_1', n_2', n_3', n_4') \le  \max_{j=1,2}\{(n_1+n_4- j)(n_2+n_3-(3- j))\}+ (k-2)n_1.
\]
If $n_1'> n_2'+n_3'$, then $g_{k-1}(n_1', n_2', n_3', n_4') = n_1'(n_2'+n_3' + n'_4)+ (k-2)(n_2' + n'_3)$. In this case, we must have $n_1=n_2+n_3-t$ for $t=0,1$, $n'_2= n_2 -1$, and $n'_3 = n_3 -1$. Thus $n'_i = n_i -1$ either for $i\in [3]$ or for $i\in \{2, 3, 4\}$, and consequently 
\begin{align*}
g_{k-1}(n_1', n_2', n_3', n_4')\le \max\{ & (n_1 -1) (n_2 + n_3 + n_4 - 2) + (k-2) (n_2 + n_3 - 2), \\
  & n_1 (n_2 + n_3 + n_4 - 3) + (k-2) (n_2 + n_3 - 2).
\end{align*}
Since $n_1=n_2+n_3-t$ for $t=0,1$, it follows that 
\[
g_{k-1}(n_1', n_2', n_3', n_4') \le  \max_{j=1,2,3}\{(n_2+n_3- (3-j))(n_1+n_4 -j)\}+ (k-2)(n_1-1).
\]
%, which further implies that $n_1\ge n_2+n_3-1$.
%Recall that $n_1\le n_2+n_3$ and note that $n_1'\le n_2'+n_3'+2$
Putting all cases together with $e(G\setminus S)\le g_{k-1}(n_1', n_2', n_3', n_4')$, we conclude that 
\begin{equation}
\label{eq:eGS}
e(G\setminus S) \le \max_{j=1,2,3}\{(n_1+n_4- j)(n_2+n_3-(3- j))\}+ (k-2)n_1.
\end{equation}

%and moreover, if the maximum is achieved for $i=3$, then $n_1=n_2+n_3-t$ for $t=0,1$.
%In fact, when $n_1'\le n_2'+n_3'$ the claim holds and the maximum can be taken over $i=1,2$ only.
%Otherwise $n_1'> n_2'+n_3'$, and we must have $n_2+n_3-1\le n_1\le n_2+n_3$ (if $n_1\le n_2+n_3-2$ then $n_1'\le n_2'+n_3'$).
%Then we can write $n_1=n_2+n_3-t$ for $t=0,1$.
%%Plug in this to the two following upper bounds of $g_{k-1}(n_1', n_2', n_3', n_4')$: 
%Note that $V_1\setminus S$ is always the largest part and thus
%\[
%g_{k-1}(n_1', n_2', n_3', n_4')\le \max\{(n_1-1)(n_2+n_3+n_4-2), n_1(n_2+n_3+n_4-3)\}+(k-2)(n_2+n_3-1)
%\]
%Now the claim follows after rewriting $n_1$ and $n_2+n_3$ for each other by $n_1=n_2+n_3-t$.
%\begin{itemize}
%%\item $e(G\setminus S) \le \max\{(n_1+n_4-1)(n_2+n_3-2), (n_1+n_4-2)(n_2+n_3-1)\}+ (k-2)n_1$, %(if $S\cap V_1=\emptyset$)
%%\item $e(G\setminus S) \le (n_1+n_4-2)(n_2+n_3-1)+ (k-2)n_1$, 
%\item $(n_1-1)(n_2+n_3+n_4-2)+(k-2)(n_2+n_3-1)$, or 
%\item $n_1(n_2+n_3+n_4-3)+(k-2)(n_2+n_3-1)$,
%\end{itemize}

%depending on the extremal number and the part sizes for $G\setminus S$ (three of them are $n_i-1$ and the other is $n_i$, for $i\in [4]$).
%In particular, when the last two cases happen, we must have $n_1\ge n_2+n_3-1$ (otherwise $n_1\le n_2+n_3-2$ and thus $n_1'\le n_2'+n_3'$).

Recall that  $e(S; G) := e(G) - e(G\setminus S)$.  We next claim that $e(S; G)\ge \frac32 \sum_{i\in [4]}n_i + 3k$.
%\begin{align}\label{eq:eSG}
%e(S; G)\ge (3/2) \sum_{i\in [4]}n_i + 3k.
%\end{align}
Indeed, if the maximum in~\eqref{eq:eGS} is achieved by $j=1,2$, then, together with~\eqref{eq:contra}, it gives
%Combining~\eqref{eq:contra} and~\eqref{eq:eGS}, we obtain
\[
%\label{eq:TG}
e(S; G)> \sum_{i\in [4]} n_i +  \min\{n_1+n_4, n_2+n_3\} + n_1-2 \ge \frac32 \sum_{i\in [4]} n_i + n_4-2\ge \frac32 \sum_{i\in [4]} n_i + 3k,
\]
where we used $n_4\ge 6k^2$ in the last inequality.
%For the other two cases, because $n_2+n_3-1\le n_1\le n_2+n_3$, we can rewrite the upper bounds as
%Note that we can rewrite
%\[
%(n_1-1)(n_2+n_3+n_4-2)+(k-2)(n_2+n_3-2) = (n_2+n_3-1)(n_1+n_4-2)+(k-2)(n_1-2).
%\]
%So we get the upper bound as in the second case, and thus~\eqref{eq:TG} holds.
%For the last case, we must have $n_2+n_3-1\le n_1\le n_2+n_3$.
%If $n_1 = n_2+n_3-1$, then we similarly rewrite 
%\[
%n_1(n_2+n_3+n_4-3)+(k-2)(n_2+n_3-2) = (n_2+n_3-1)(n_1+n_4-2)+(k-2)(n_1-1)
%\] 
%and get the upper bound as in the second case, and thus~\eqref{eq:TG} holds.
Otherwise, the maximum in~\eqref{eq:eGS} is achieved by $j=3$, that is, $e(G\setminus S) \le (n_1+n_4-3)(n_2+n_3)+ (k-2)n_1$.
% and we have $n_2+n_3\ge n_1$, and we infer $(5/2)n_1 \ge n_1+n_2+n_3+n_4 $.
%Thus, $(15/4)n_1\ge (3/2)\sum_{i\in [4]}n_i $.
%So we obtain $e(G\setminus S) \le n_1(n_2+n_3+n_4-3)+(k-2)(n_1-1)$.
By~\eqref{eq:contra}, we get
\begin{align*}
e(S; G)&> (n_1+n_4)(n_2+n_3)+ (k-1)n_1 - (n_1+n_4-3)(n_2+n_3)- (k-2)n_1 \\
&= n_1+3n_2 +3n_3 \ge \frac32 \sum_{i\in [4]}n_i + \frac{n_4}2 \ge \frac32 \sum_{i\in [4]}n_i + 3k,
\end{align*}
where we used the assumption $n_2+n_3\ge n_1$ and $n_2, n_3\ge n_4$.
%Since $n_4\ge 6k^2$, $e(S; G)\ge (3/2) \sum_{i\in [4]}n_i + 3k$ holds for all cases.

Let $S=xyz$ and note that $d(x)+d(y)+d(z) = e(S; G)+3$.
By averaging, without loss of generality, we may assume that
\[
d(x)+d(y) \ge \frac23 \left(\frac32 \sum_{i\in [4]} n_i + 3k \right) = \sum_{i\in [4]} n_i + 2k.
\]
By the moreover part of Fact~\ref{fact:com_neigh}, $xy$ is rich and we are done.
\end{proof}

For two disjoint sets $A, B\subseteq V(G)$, let $d(A, B) = e(A, B)/(|A| |B|)$ be the density of the bipartite graph with parts $A$ and $B$.
A pair $(V_i, V_j)$ is called \emph{full} if $d(V_i\setminus Z, V_j)=d(V_j\setminus Z, V_i)=1$; $(V_i, V_j)$ is called \emph{empty} if $e(V_i\setminus Z, V_j)=e(V_i, V_j\setminus Z)=0$. We have the following observation.

\begin{obs}\label{obs:empty}
%Suppose~\eqref{eq:contra} holds.
For distinct $i, j, t\in [4]$, if $d(V_i\setminus Z, V_j)=d(V_i\setminus Z, V_t)=1$, then $(V_j, V_t)$ must be empty because any edge in $(V_j,V_t)$ but not in $(V_j\cap Z, V_t\cap Z)$ will create a triangle with at most one vertex in $Z$, contradicting~\ref{item:3}.
In particular, if both $(V_i, V_j)$ and $(V_i, V_t)$ are full, then $(V_j, V_t)$ is empty.
\end{obs}

\begin{clm}
\label{claim:1} 
%Suppose~\eqref{eq:contra} holds.
Fix $i\ne j\in [4]$.
If $d(x)+d(y) \ge \sum_{i\in [4]}n_i$ for every edge $x y\in V_i\times V_j$, 
%for $\{i,j\}\in \{\{1,2\}, \{1,3\},\{2,4\}, \{3,4\}\}$.
then either %one of the following holds,
\begin{itemize} 
\item $e(V_i\setminus Z, V_j\setminus Z)=0$ (this is weaker than $(V_i, V_j)$ being empty) or 
\item $d(V_i\setminus Z, V_j)=d(V_j\setminus Z, V_i)=1$, and $d(x)+d(y) = \sum_{i\in [4]}n_i$.
\end{itemize}
Moreover, if $d(x)+d(y) > \sum_{i\in [4]}n_i$ for every edge $x y\in V_i\times V_j$, then $(V_i, V_j)$ is empty.
\end{clm}

\begin{proof}
%Without loss of generality, let $i=1$ and $j=2$.
Assume that $\{i,j,t,\ell\}=[4]$.
Suppose there is an edge $xy \in (V_i\setminus Z)\times (V_j\setminus Z)$.
Note that if $x$ and $y$ have a common neighbor $z$, then as $x, y\notin Z$, none of the edges of $xyz$ is rich, contradicting~\ref{item:3}.
Thus, $x$ and $y$ have no common neighbor.
By Fact~\ref{fact:com_neigh}, $d(x)+d(y)\le \sum_{i\in [4]} n_i$.
If $d(x)+d(y) \ge \sum_{i\in [4]}n_i$, then Fact~\ref{fact:com_neigh} implies that $V_j \subseteq N(x)$ and $V_i \subseteq N(y)$.
In particular, $x y'\in E(G)$ for every $y'\in V_j\setminus Z$. Applying the same argument to the edge $x y'$, we obtain that $V_i \subseteq N(y')$.
Similarly, we can derive that $V_j \subseteq N(x')$ for every $x'\in V_i\setminus Z$. Thus, %Repeated applications of this observation on these edges imply that for any $xy\in V_i\times V_j$ such that at least one of them is not in $Z$, $xy\in E(G)$, that is, 
$d(V_i\setminus Z, V_j)=d(V_j\setminus Z, V_i)=1$.

Now assume $d(x)+d(y) > \sum_{i\in [4]}n_i$ for every edge $x y\in V_i\times V_j$.
If $e(V_i\setminus Z, V_j\setminus Z) \ne 0$, then the arguments in the previous paragraph provide a contradiction.
Suppose there is an edge $xy\in (V_i\cap Z)\times (V_j\setminus Z)$.
As $d(x)+d(y) > \sum_{i\in [4]}n_i$, $x$ and $y$ have some common neighbors in $V_t\cup V_\ell$.
But since $y\notin Z$, by~\ref{item:3}, their common neighbors must be in $(V_t\cup V_\ell)\cap Z$.
Since $e(V_i\setminus Z, V_j\setminus Z)=0$, we know that $N(y)\cap V_i\subseteq V_i\cap Z$.
Altogether, we obtain that $d(x)+d(y)\le n_j+n_t+n_\ell+|Z| < \sum_{i\in [4]}n_i$, a contradiction.
Analogous arguments show that there is no edge in $(V_i\setminus Z)\times (V_j\cap Z)$.
Thus, $e(V_i\setminus Z, V_j)=e(V_i, V_j\setminus Z)=0$, that is, $(V_i, V_j)$ is empty.
\end{proof}

%Finally we make some preparations about the progressive induction.
Consider a set $T\subseteq V(G)$ defined in \textbf{Cases~1--4} and let $n'_1, n'_2, n'_3, n'_4$ denote the sizes of the parts of $G\setminus T$.
%Then $G\setminus T\in \fA_{n'_1+\cdots+ n'_4}$.
Then $\phi(G) < \phi(G\setminus T)$ is equivalent to 
\[
e(G) -  g_k(n_1, n_2, n_3, n_4) < e(G\setminus T) - g_k(n'_1, n'_2, n'_3, n'_4),
\]
or $e(T; G)< g_k(n_1, n_2, n_3, n_4) - g_k(n'_1, n'_2, n'_3, n'_4)$.
We will prove by contradiction, assuming that $\phi(G) \ge \phi(G\setminus T)$, equivalently, 
%such that \emph{all parts of $G\setminus T$ have sizes at least $6k^2$},
\begin{equation}\label{eq:etg}
e(T; G) \ge (n_1+n_4)(n_2+n_3) + (k-1)n_1 - g_k(n'_1, n'_2, n'_3, n'_4)
\end{equation}
for every $T\subseteq V(G)$ defined in \textbf{Cases~1--4}.

%Recall that we need to prove~\ref{item:bbb} for different types of $T$ as specified in the cases.
The case when $T=\{v\}\subseteq V_1$ occurs in all four cases so we consider it before the cases.
Since $n_1\le n_2+n_3$, we have three possibilities:
\begin{itemize}
\item if $n_1>n_2$, then $g_k(n_1-1, n_2, n_3, n_4) = (n_1-1+n_4)(n_2+n_3) + (k-1)(n_1-1)$;
\item if $n_1=n_2>n_4$, then $g_k(n_1-1, n_2, n_3, n_4) = (n_1+n_4)(n_2+n_3-1) + (k-1)n_1$;
\item if $n_1=n_4$, then $g_k(n_1-1, n_2, n_3, n_4) = (n_1+n_4-1)(n_2+n_3) + (k-1)n_1$;
\end{itemize}
Thus~\eqref{eq:etg} implies that for every $v\in V_1$,
% and $e(\{v\}; G) \ge (n_1+n_4)(n_2+n_3) + (k-1)n_1 - g_k(n_1-1, n_2, n_3, n_4)$, then
%  \begin{enumerate}[label=(A\arabic*)]
%  \item $d(v) \ge n_2 + n_3+k-1$, if $n_1 > n_2$; \label{item:A1}
%  \item $d(v) \ge n_1 + n_4$, if $n_1 = n_2$. \label{item:A2}
%  \end{enumerate}
  \begin{align}
  \label{eq:dV1}
  d(v)\ge 
  \left\{\begin{array}{lr}
         n_2 + n_3+k-1, & \text{if } n_1 > n_2, \\
        n_1 + n_4,  & \text{if } n_1 = n_2.
        \end{array}\right.
  \end{align}
  
%\item for any $v\in V_i$, $i=2,3$, \label{item:B}
%  \begin{enumerate}[label=(B\arabic*)]
%  \item $d(v) \ge n_1 + n_4$, if $n_i > n_4$ and $n_1<n_2+n_3$; \label{item:B1}
%  \item $d(v) \ge n_2 + n_3$, otherwise. \label{item:B2}
%  \end{enumerate}
%\item (Only for \textbf{Case 2}) for any $v\in V_4$, $d(v) \ge n_2 + n_3$.  \label{item:C}
%\end{enumerate}

%To illustrate the cases for $G\setminus \{v\}$, 
%\item $g_k(n_1, n_2, n_3, n_4-1) = (n_1+n_4-1)(n_2+n_3) + (k-1)n_1$.
%\end{enumerate}
%Next, for $i=2,3$,
%\begin{enumerate}[resume]
%\item if $n_1< n_2+n_3$ and $n_i>n_4$, then $g_k(n_1, n_i-1, n_{5-i}, n_4) = (n_1+n_4)(n_2+n_3-1) + (k-1)n_1$;
%\item if $n_i=n_4$, then $g_k(n_1, n_i-1, n_{5-i}, n_4) = (n_1+n_4-1)(n_2+n_3) + (k-1)n_1$;
%\item if $n_1= n_2+n_3$ and $n_i>n_4$, then $g_k(n_1, n_i-1, n_{5-i}, n_4) = n_1(n_2+n_3+n_4-1) + (k-1)(n_2+n_3-1)$.

%Then~\ref{item:A1} and~\ref{item:A2} follow from these equations and~\eqref{eq:etg}.
%In other words, if~\ref{item:A1} or~\ref{item:A2} fails in the corresponding case, then~\ref{item:bbb} holds for $T=\{v\}\subseteq V_1$ and we are done.
%By \eqref{eq:etg'} and straightforward calculations we see that (1) implies~\ref{item:A1}, (2) and (3) imply~\ref{item:A2} and (4) implies~\ref{item:C}.
%Similarly, (5) implies~\ref{item:B1} and (6) and (7) imply~\ref{item:B2}.

\subsection{Proof of Cases 1--4}
\label{sec:22}

After these preparations, we return to the proof of \textbf{Cases 1--4}.
%We restate the cases for convenience.
Recall that $n_1\le n_2+n_3$ in all these cases.
Recall also that $n_i\ge 6k^2$ for $i\in [4]$, so we can always assume that $V_i\setminus Z\neq\emptyset$.
Moreover, by~\eqref{eq:M0}, we have $M_0(k)\ge N_0(k-1)+3$, and thus we can apply the induction hypothesis \ref{item:IH} on any \emph{$(k-1)K_3$-free} subgraph $G\setminus S$, whenever $|S|\le 3$ (and thus $v_4(G\setminus S)\ge 6k^2 - 3 \ge 6(k-1)^2$).

\medskip
\noindent\textbf{Case 1. $n_1 > n_3$ and $n_2> n_4$.} 

%\begin{proof}%[Proof of {\rm \textbf{Case 1}}]

In this case \eqref{eq:etg} holds for every crossing set $T= xy\in V_1\times V_2$. %of size at most two.
%For two vertices $x\in V_1$, $y\in V_2$, 
Since the part sizes of $G\setminus \{x,y\}$ are $n_1-1\ge \{n_2-1, n_3\} \ge n_4$.
%So either $\phi(G)< \phi(G\setminus \{x,y\})$ (then we are done by choosing $T=xy\in V_1\times V_2$), or we have
By~\eqref{eq:etg}, we have
\begin{align*}
e(xy; G) &\ge (n_1+n_4)(n_2+n_3) + (k-1)n_1 - ((n_1+n_4-1)(n_2+n_3-1) + (k-1)(n_1-1)) \\
&= \sum_{i\in [4]} n_i + k-2.
\end{align*}
If $xy\in E(G)$, then $d(x)+d(y) =  e(xy; G) + 1\ge \sum_{i\in [4]} n_i + k-1>\sum_{i\in [4]} n_i$.
%Then either $\phi(G)< \phi(G\setminus \{x,y\})$ and we are done, or we have $d(x)+d(y) >\sum n_i$ for every edge $x y\in V_1\times V_2$.
By Claim~\ref{claim:1}, $(V_1, V_2)$ is empty.
For every $x\in V_1\setminus Z$, we thus have $d(x)\le n_3+n_4<\min\{n_2+n_3, n_1+n_4\}$, contradicting~\eqref{eq:dV1}.
%\end{proof}

\medskip
\noindent\textbf{Case 2. $n_1 = n_2 = n_3 \ge n_4 > 6k^2$.} 

%\begin{proof}%[Proof of {\rm \textbf{Case 2}}]
In this case \eqref{eq:etg} holds for any one-element set $T\subset V(G)$.
%For two vertices $x\in V_1\cup V_2\cup V_3$, $y\in V_4$, the sizes of four parts of $G\setminus \{x,y\}$ are $n_1\ge \{n_2-1, n_3\} \ge n_4-1$.
%%So either $\phi(G)< \phi(G\setminus \{x,y\})$ (then we are done by choosing $T=xy$), or we have
%By~\eqref{eq:etg}, %we have
%\begin{align*}
%e(xy; G) &\ge (n_1+n_4)(n_2+n_3) + (k-1)n_1 - ((n_1+n_4-1)(n_2+n_3-1) + (k-1)n_1) \\
%&= \sum_{i\in [4]} n_i -1.
%\end{align*}
%In addition, if $xy\in E(G)$, then we have $d(x)+d(y) \ge \sum_{i\in [4]} n_i$.
Write $n_1=n_2=n_3=n$. %and note that $n_1<n_2+n_3$.
For any $x\in V_1\cup V_2\cup V_3$, by~\eqref{eq:etg}, we have
\[
d(x)=e(\{x\}; G) \ge 2n(n+n_4) + (k-1)n - g_k(n, n, n-1, n_4), 
\]
where $g_k(n, n, n-1, n_4) = (2n-1)(n+n_4)+ (k-1)n$ if $n>n_4$ and $g_k(n, n, n-1, n_4) = 2n(n+n_4-1)+ (k-1)n$ if $n=n_4$.
Thus, we have $d(x)\ge \min\{2n, n+n_4\}=n+n_4$.
Similarly, for $y\in V_4$, by~\eqref{eq:etg}, we have
\begin{align}\label{eq:dy}
d(y)=e(\{y\}; G) \ge 2n(n+n_4) + (k-1)n - \big(2n(n+n_4-1) + (k-1)n \big) = 2n.
\end{align}
These together imply $d(x)+d(y) \ge \sum n_i$ for every edge $x y\in (V_1\cup V_2\cup V_3)\times V_4$.
For $i=1, 2, 3$, Claim~\ref{claim:1} implies that either $(V_i, V_4)$ is full or $e(V_i\setminus Z, V_4\setminus Z)=0$. 
If $e(V_i\setminus Z, V_4\setminus Z)=0$ holds for at least two values of $i\in \{1, 2, 3\}$, then for every $y\in V_4\setminus Z$, we have $d(y)\le n + |Z| < 2n$ (as $n\ge M_0(k) /4 > 6 k^2$), contradicting \eqref{eq:dy}.

This implies that at least two of $(V_1, V_4)$, $(V_2, V_4)$, and $(V_3, V_4)$ must be full.
Without loss of generality, assume $(V_1, V_4)$ and $(V_2, V_4)$ are full.
By Observation~\ref{obs:empty}, $(V_1, V_2)$ is empty.
Next, we claim that $(V_3, V_4)$ is empty.
Indeed, let $x\in V_2\setminus Z$ and recall that $d(x)\ge n+n_4$.
%let $x\in V_2\setminus Z$ and by~\ref{item:B}, we have $d(x)\ge \min\{2n, n+n_4\} = n+n_4$.
Since $(V_1, V_2)$ is empty, we have $d(x)\le n+n_4$.
Thus, $d(x)=n+n_4$ and in particular $V_3\subseteq N(x)$. %$d(x, V_3)=n$.
Since this holds for every $x\in V_2\setminus Z$, it follows that $d(V_2\setminus Z, V_3)=1$.
Thus $(V_3, V_4)$ is empty by Observation~\ref{obs:empty}.
Together with~\ref{item:2}, we infer 
\[
e(G)= e(G[Z]) + e(V\setminus Z; G) < \binom{|Z|}2 + (n_1+n_2)(n_3+n_4)\le (n_1+n_2)(n_3+n_4) + (k-1)n_1, 
\]
contradicting~\eqref{eq:contra},
The previous inequality follows from $ \binom{|Z|}2 \le 18(k-1)^4\le (k-1)n_1$, which follows from $n_1\ge M_0(k)/4$ and~\eqref{eq:M0}.
%\end{proof}

%The remaining two cases are more involved.
\medskip
\noindent\textbf{Case 3. $n_1 = n_2 = n_3 > n_4 = 6k^2$.} 

%\begin{proof}%[Proof of {\rm \textbf{Case 3}}]
%We may apply \eqref{ a set $T\subset V_1\cup V_2\cup V_3$.
Write $n_1=n_2=n_3=n$. We assume that 
\begin{align}\label{eq:n1}
n_1 \ge 30k^2,
\end{align}
as otherwise $\sum n_i\le 3\cdot 30k^2 + 6k^2 \le M_0(k)$ by \eqref{eq:M0}, contradicting the assumption $|G|> M_0(k)$.
By~\eqref{eq:dV1} and the similarity of $V_1, V_2$, and $V_3$, we have $d(x) \ge n+n_4$ for every $x \in V_1\cup V_2\cup V_3$.
We claim that for $y\in V_4$,
\begin{equation}
\label{eq:dyV4}
d(y)\le 2n+2k-1.
\end{equation}
Otherwise, pick $k$ neighbors $x_1, \dots, x_{k}$ of $y$ from the same part of $G$.
For each $i$, since $d(x_i) \ge n+n_4$, we have $d(x_i)+d(y) \ge \sum n_i+2k-1$, yielding that $x_iy$ is rich by Fact~\ref{fact:com_neigh}.
However, this contradicts~\ref{item:1}.
 
%Our first goal is to show that if $G[V_1\cup V_2\cup V_3]$ contains a triangle, then we can apply the standard induction.

\medskip
\noindent \textbf{Claim.} The graph $G[V_1\cup V_2\cup V_3]$ is $K_3$-free.
%\emph{Every triangle $xyz\in V_1\times V_2\times V_3$ contains at least two rich edges}. 

%\emph{Proof of the claim.} 
\begin{proof} %[Proof of the claim]
Suppose instead, there exists a triangle $xyz\in V_1\times V_2\times V_3$. Without loss of generality, assume that $d(x)\ge d(y)\ge d(z)$.
We first claim that
\begin{align}\label{eq:xyz}
d(x)+d(y)+d(z) \ge 5n+2n_4+k.
\end{align}
Otherwise $d(x)+d(y)+d(z)\le 5n+2n_4+k-1$ and $e(xyz; G)= d(x)+d(y)+d(z) - 3 \le 5n+2n_4+k-4$.
Then, by~\eqref{eq:contra}, 
\begin{align*}
e(G\setminus \{x, y, z\}) &= e(G) - e(xyz; G) > g_k(n, \,n,n,n_4) - (5n+2n_4+k-4)\\
&= 2n(n+n_4) +(k-1)n - (5n+2n_4+k-4)\\
&=  (2n-2)(n-1+n_4) + (k-2)(n-1) \\
&=g_{k-1}(n-1, n-1, n-1, n_4).
\end{align*}
%Aiming to apply the induction on $G-xyz$, we compute
%\begin{align*}
%g_k(&n, \,n,n,n_4) - g_{k-1}(n-1, n-1, n-1, n_4) \\
%&= 2n(n+n_4) +(k-1)n - ((2n-2)(n-1+n_4) + (k-2)(n-1))\\
%& = 5n+2n_4+k-4.
%\end{align*}
By induction hypothesis \ref{item:IH}, we obtain a copy of $(k-1) K_3$ in $ G\setminus \{x, y, z\} $. Together with the triangle $xyz$, this contradicts the assumption $G$ is $k K_3$-free.

 %for every triangle $xyz\in V_1\times V_2\times V_3$.

We next claim that at least two of $xy, yz, xz$ are rich and thus all $x, y, z\in Z$.
%every triangle in $G[V_1\cup V_2\cup V_3]$ contains \emph{two} rich edges.
%Suppose there is such a triangle $xyz$ and 
Indeed, if $d(x)<2n+n_4-k$, then by \eqref{eq:xyz}, 
\[
d(y)+d(z)> 5n+2n_4+k- (2n+n_4-k) = 3n+n_4+2k > \sum n_i + 2k-1.
\]
By Fact~\ref{fact:com_neigh}, $yz$ is rich.
Since $d(x)$ is the largest, this argument implies that all three edges of $xyz$ are rich, as desired.
Otherwise, $d(x)\ge 2n+n_4-k$ and recall that $d(y)\ge d(z)\ge n+n_4$.
Thus 
\[
d(x)+d(y) \ge d(x)+d(z)\ge 3n+2n_4-k \ge \sum n_i + 2k-1
\]
because $n_4= 6k^2\ge 3k-1$. By Fact~\ref{fact:com_neigh}, both $xy$ and $xz$ are rich, as desired.

%Assume that $xyz\in V_1\times V_2\times V_3$ is a triangle with $d(x)\ge d(y)\ge d(z)$, 
The claim in the previous paragraph applies to all triangles in $V_1\cup V_2 \cup V_3$. Therefore, all the common neighbors of $x$ and $y$ in $V_1\cup V_2 \cup V_3$ are in $Z$ and consequently, $| N(x)\cap N(y)|\le |Z|+|V_4|\le 6k^2+n_4$, and consequently, $d(x)+d(y)\le \sum n_i + 6k^2+n_4 = 3n+2n_4+ 6k^2$.
%In particular, the fact above implies that if $xyz$ is a triangle in $G[V_1\cup V_2\cup V_3]$, then either $xz$ or $yz$ must be rich, that is, $z\in N_R(x)\cup N_R(y)$.
%So we infer that for every edge $xy$ in e.g. $V_1\cup V_2$, $N_G(x)\cap N_G(y)\cap V_3\subseteq (N_R(x)\cup N_R(y))\cap V_3$.
%Since by \ref{item:1} both $x$ and $y$ send at most $k-1$ rich edges to $V_3$, $x$ and $y$ have at most $2(k-1)$ common neighbors $z$ in $V_3$, implying that $d(x)+d(y)\le 3n+2n_4+2(k-1)$.
On the other hand, \eqref{eq:xyz} and the assumption $d(x)\ge d(y)\ge d(z)$ imply that 
\begin{equation}\label{eq:xyz-rich}
d(x)+d(y)\ge \frac23 (5n+2n_4+k) = \frac{10}3n + \frac43 n_4 + \frac23k > 3n+2n_4+ 6k^2
\end{equation}
because $n\ge 30k^2 = 2n_4+18k^2$ by~\eqref{eq:n1}. This gives a contradiction. 
%Together with~\eqref{eq:xyz-rich}, $xy$ can not be in a triangle in $G[V_1\cup V_2\cup V_3]$, yielding the claim.
\end{proof}

%\medskip
%Now suppose there is a triangle $xyz\in V_1\times V_2\times V_3$ with $d(x)\ge d(y)\ge d(z)$, then since $d(x)+d(y)+d(z) \ge 5n+2n_4+k$, we infer that
%\[
%d(x)+d(y)\ge \frac23 (5n+2n_4+k) = \frac{10}3n + \frac43 n_4 + \frac23k > 3n+2n_4+2k,
%\]
%as $n\ge 144k^4/(k-1)\ge 2n_4+4k$, a contradiction.
By the claim, $G[V_1\cup V_2\cup V_3]$ is $K_3$-free, and thus has at most $2n^2$ edges by Theorem~\ref{thm:BES}.
Together with~\eqref{eq:dyV4} and~\eqref{eq:n1}, we obtain that 
\[
e(G)\le 2n^2 + n_4\cdot (2n+2k-1) = 2n(n+n_4) +(2k-1)n_4 < 2n(n+n_4) +(k-1)n,
\] 
contradicting \eqref{eq:contra}.

\medskip
\noindent\textbf{Case 4. $n_1 > n_2 = n_3 = n_4 $.} %We verify~\ref{item:bbb} for $T\subset V(G)$.

%\begin{proof}%[Proof of {\rm \textbf{Case 4}}]
Assume $n_2=n_3=n_4=n$ and recall that $n_1\le 2n$.
We first claim that 
\begin{align}\label{eq:dy4}
d(x)\le 3n \text{ for all } x\in V_1,  \  \text{and} \ d(y)\le n_1+n+k-1 \text{ for all } y\in  V_2\cup V_3\cup V_4.
\end{align}
Indeed, the bound $d(x)\le 3n$ for $x\in V_1$ is trivial. Suppose to the contrary, that there is a vertex $y\in V_2\cup V_3\cup V_4$ with $d(y)\ge n_1+n+k$.
It follows that $|N(y)\cap V_1| \ge d(y) - 2n \ge k$.
Assume that $x_1,\dots, x_k\in N(y)\cap V_1$.
By~\eqref{eq:dV1}, we have $d(x_j)\ge 2n+k-1$.
Thus, we infer that $d(x_j)+d(y)\ge n_1+3n+ 2k-1$. By Fact~\ref{fact:com_neigh}, we have $x_1 y, \dots, x_k y\in E(R)$.
However, this contradicts~\ref{item:1}.
%We will use this below together with the trivial bound $d(x)\le 3n$ for any $x\in V_1$.

%Let $e\in E(R)$ be a rich edge.
%By definition, given any set $S\subseteq V(G)$ that spans a copy of $(k-1)K_3$, $e$ must intersect $S$, as otherwise we can find a triangle which contains $e$ and does not intersect S$, a contradiction.
%This implies that $G\setminus e$ is $(k-1)K_3$-free and below we seek to show $\phi(G)=0$ by applying the induction hypothesie on $G\setminus e$.

We next claim that there is no rich edge in $V_1 \times (V_2\cup V_3\cup V_4)$. Suppose to the contrary, that $xy\in V_1 \times (V_2\cup V_3\cup V_4)$ is a rich edge. By \eqref{eq:dy4}, we have $e(xy; G)=d(x)+d(y)-1\le n_1+4n+k-2$. By \eqref{eq:contra}, it follows that
\begin{align*}
e(G\setminus \{x, y\}) &= e(G) -e(xy; G) > %g_k(n_1, \,n,n,n) - (n_1+4n+k-2)\\ &=
2n (n_1+n)  +(k-1)n_1 -  (n_1+4n+k-2) \\
&= 2n (n_1+n-2)  + (k-2)(n_1-1) \\
&= g_{k-1}(n_1-1, n, n, n-1).
\end{align*}
By induction hypothesis \ref{item:IH}, $G\setminus \{x, y\}$ contains a copy $S$ of $(k-1) K_3$. Since $xy$ is rich, we can find a triangle in $G\setminus S$ containing $xy$, contradicting the assumption that $G$ is $k K_3$-free. 

%\begin{align*}
%g_k(&n_1, \,n,n,n) - g_{k-1}(n_1-1, n, n, n-1) \\
%&= (n_1+n) 2n +(k-1)n_1 - ((n_1+n-2) 2n + (k-2)(n_1-1))\\
%& = n_1+4n+k-2.
%\end{align*}
%By the upper bounds on degrees, we have 
%By induction, $e(G\setminus e)\le g_{k-1}(n_1-1, n, n, n-1)$, which yields $e(G)\le g_{k}(n_1, n, n, n)$, contradicting~\eqref{eq:contra}.
%Thus, we may assume that there is no rich edge from $V_1 \times (V_2\cup V_3\cup V_4)$.
%In particular, $A\cap V_1=\emptyset$.

Now we show that there is no triangle intersecting $V_1$.
Suppose to the contrary, there is a triangle $xyz$ with $x\in V_1$. %$xyz\in V_1\times V_2\times V_3$.
If $d(x)+d(z) \ge n_1+3n+ 2k-1$, then, by Fact~\ref{fact:com_neigh}, $xy$ is rich, contradicting our earlier claim. 
%$x$ and $z$ have at least $2k-1$ common neighbors, and thus at least $k$ common neighbors in the same part of $G$. This implies that $xz\in E(R)$, 
We thus assume that $d(x)+d(z) < n_1+3n+2k-1$. Together with \eqref{eq:dy4}, it gives that $d(x)+d(y)+d(z) < 2n_1+4n+3k-2$, and
$e(xyz; G)= d(x)+d(y)+d(z)-3 < 2n_1+ 4n+ 3k- 5$.
By \eqref{eq:contra}, it follows that
\begin{align*}
e(G\setminus \{x, y, z\}) &= e(G) - e(xyz; G) %> g_k(n_1, \,n,n,n) - (2n_1+ 4n+ 3k- 8)\\
> 2n (n_1+n) +(k-1)n_1 -  (2n_1+ 4n+ 3k- 5) \\
&= (n_1+n-2)(2n-1) + (k-2)(n_1-1) + n - 2k + 1\\
&= g_{k-1}(n_1-1, n, n-1, n-1) + n - 2k + 1.
\end{align*}
%because $g_{k-1}(n_1-1, n, n-1, n-1) = (n_1+n-2)(2n-1) + (k-2)(n_1-1)$. 
By \ref{item:IH}, $G\setminus \{x, y, z\}$ contains a copy of $(k-1) K_3$. Together with the triangle $xyz$, this contradicts the assumption that $G$ is $k K_3$-free.

%Suppose there is a triangle $xyz\in V_1\times V_2\times V_3$.
%Note that
%\begin{align*}
%g_k(&n_1,\,n,n,n) - g_{k-1}(n_1-1, n, n-1, n-1) \\
%&= (n_1+n) 2n +(k-1)n_1 - ((n_1+n-2)(2n-1) + (k-2)(n_1-1)) \\
%&= 2n_1+5n+k-4.
%\end{align*}
%Thus, if $e(xyz; G) \le 2n_1+5n+k-4$, then similarly by $e(G\setminus \{x,y,z\})\le g_{k-1}(n_1-1, n, n-1, n-1)$ we derive a contradiction with~\eqref{eq:contra}.
%So we may assume that $d(x)+d(y)+d(z)=e(xyz; G)+3 \ge 2n_1+5n+k$.
%Since $d(y)\le n_1+n+k-1$, we obtain that 

We assumed that $G$ contains $k-1$ disjoint triangles. Let $T_1$ be a triangle of $G$.
%is maximal and $k\ge 2$, $G$ contains a triangle.
By the claim of the previous paragraph, $T_1$ must be in $V_2\cup V_3\cup V_4$.
%However, we will show that this is also not possible, which finishes the proof of Case 3.
Moreover, by~\ref{item:3}, $T_1$ must contain a rich edge $xy$.
%Without loss of generality, suppose $xy\in V_2\times V_3$.
Below we show that 
\begin{align}\label{eq:Gxy}
e(G\setminus \{x, y\}) > g_{k-1}(n_1, n, n-1, n-1). 
\end{align}
Then, by \ref{item:IH}, $G\setminus \{x, y\}$ contains a copy $S$ of $(k-1) K_3$. Since $xy$ is rich, we can find a triangle in $G\setminus S$ containing $xy$, contradicting the assumption that $G$ is $k K_3$-free. 

We first assume that $n_1=2n$.
If $d(x)+d(y) > 6n$, then $x$ and $y$ have a common neighbor in $V_1$, contradicting the earlier claim that there is no triangle intersecting $V_1$.
We thus assume that $d(x) + d(y)\le 6n$. Thus $e(xy; G) \le 6n-1$. By \eqref{eq:contra}, it follows that
\begin{align*}
e(G\setminus \{x, y\})  &> g_k(2n, n,n, n) - (6n-1) \\
&=  3n\cdot 2n + 2n (k-1) -  (6n-1) \\
&=2n(3n-2) + (k-2)(2n-1) + k-1 \\
&= g_{k-1}(2n, n, n-1, n-1) + k-1.
\end{align*}
Thus \eqref{eq:Gxy} holds. 
%because $g_{k-1}(2n, n, n-1, n-1) = 2n(3n-2) + (k-2)(2n-1)$.
%By \ref{item:IH}, $G\setminus \{x, y\}$ contains a copy $S$ of $(k-1) K_3$. Since $xy$ is rich, we can find a triangle in $G\setminus S$ containing $xy$, contradicting the assumption that $G$ is $k K_3$-free. 
Second, assume $n_1<2n$. %By \eqref{eq:dy4}, we have $d(x)+d(y)\le 2(n_1+n+k-1) < 2n_1+3n+k-3$.
By \eqref{eq:dy4}, we have $e(xy; G) = d(x)+d(y)-1 \le 2(n_1+n+k-1) -1$. By \eqref{eq:contra}, it follows that
\begin{align*}
e(G\setminus \{x, y\}) &> g_k(n_1, \,n,n,n) - (2n_1+ 2n+ 2k-3)\\
&=  (n_1+n) 2n +(k-1)n_1 -  (2n_1+ 2n+ 2k-3)\\
&= (n_1+n-1)(2n-1) + (k-2) n_1 + n -2k + 2\\
&= g_{k-1}(n_1, n, n-1, n-1) + n- 2k+2.
\end{align*}
Thus \eqref{eq:Gxy} holds.
%By \ref{item:IH}, $G\setminus \{x, y\}$ contains a copy $S$ of $(k-1) K_3$. Since $xy$ is rich, we can find a triangle in $G\setminus S$ containing $xy$, contradicting the assumption that $G$ is $k K_3$-free. 
%
%Note that
%\begin{align*}
%g_k(&n_1,\,n,n,n) - g_{k-1}(n_1, n, n-1, n-1) \\
%&= (n_1+n) 2n +(k-1)n_1 - ((n_1+n-1)(2n-1) + (k-2)(n_1-1)) \\
%&= 2n_1+3n+k-3.
%\end{align*}
%
%Since by induction we have $e(G\setminus \{x, y\})\le g_{k-1}(n_1, n, n-1, n-1) $, we derive a contradiction with~\eqref{eq:contra}.

The proof of Theorem~\ref{thm:key} is now completed.
\end{proof}

\section{Concluding remarks}

In this paper we solved Problem~\ref{pro1} for $r=4$ and $t=3$ when all $n_i$'s are large.
%The idea in our proof should be helpful for other cases of Problem~\ref{pro1}.
The idea in our proof should be helpful for proving Conjecture~\ref{conj} in general.
However, to determine the maximum in \eqref{eq:conj}, there are quite a few cases to consider even when $r=5$ and $t=3$.
Indeed, suppose $n_1\ge n_2\ge \cdots \ge n_5$ and $\{I, I'\}$ is the bipartition of $[5]$ that attained the maximum in \eqref{eq:conj}. Assume $1\in I$. Depending on the values of $n_1, \dots, n_5$, it is possible to have 
\[
I = \{1\}  \text{ or } \{1, 2\}  \text{ or } \{1, 3\}  \text{ or } \{1, 4\}  \text{ or } \{1, 5\} \text{ or } \{1, 4, 5\}.
\] 
%, so new ideas are expected.
%Indeed, denote the parts of a 5-partite graph $G$ by $V_1, \dots, V_5$ and the bipartition  by $(X, V(G)\setminus X)$. Suppose that $
%Then the possible choices for $X$ (up to symmetry) are $V_1$, $V_1\cup V_i$, $i=2,3,4,5$, and $V_1\cup V_4\cup V_5$. 
%It looks hard to analyze the inductive step for these cases (note that we also need to consider the boundary cases).

Another open problem is to find the smallest $N_0(k)$ such that Theorem~\ref{thm:main} holds.
The $N_0(k)$ provided in our proof is a double exponential function of $k$. Indeed, by \eqref{eq:M0} and $N_0(1)=1$, we have $M_0(2)= 96\cdot 2^2 = 384$ and $N_0(2)= 384^2$. It is easy to see that $N_0(k) = ( N_0(k-1) + 3)^2$ for $k\ge 3$. Thus $N_0(k-1)^2 \le N_0(k)\le 2 N_0(k-1)^2$ for $k\ge 3$.
It follows that 
\[
N_0(2)^{2^{k-2}} \le N_0(k) \le \big(2 N_0(2)\big)^{2^{k-2}}.
\]
It is interesting to know whether one can reduce $N_0(k)$ to a polynomial function (or even a linear function) of $k$.

%At last, it is also natural to consider the extension of Theorem~\ref{thm:mindeg} to other cliques $K_t$, $t\ge 4$. 
%
%\begin{conj}\label{conj1}
%Let $G$ be an $r$-partite graph whose parts have sizes $n_1,\dots, n_r$. 
%If
%\[
%\delta(G)> \max_{\mathcal P} \min_{I\in \mathcal P} n_{[r]\setminus I},
%\]
%where the maximum is over all partitions $\mathcal P$ of $[r]$ into $t-1$ parts, then $G$ contains a copy of $K_t$.
%\end{conj}

\medskip
\noindent \textbf{Acknowledgements.} We would like to thank Chunqiu Fang and Longtu Yuan for valuable feedbacks on an earlier version of the manuscript and thank Ming Chen, Jie Hu and Donglei Yang for helpful discussions.
We also thank two anonymous referees for their helpful comments that improved the presentation of this paper.

\bibliographystyle{abbrv}
\bibliography{Sep2016}

\noindent
\end{document}